\documentclass{amsart}
\usepackage{amsmath}
\usepackage{amsfonts}
\usepackage{graphics}
\usepackage{epsfig}
\usepackage{amssymb}
\usepackage{amscd}

\newtheorem{thm}{Theorem}[section]

\newtheorem{lem}[thm]{Lemma}
\newtheorem{prop}[thm]{Proposition}
\theoremstyle{definition}
\newtheorem{Def}[thm]{Definition}
\newtheorem*{ack}{Acknowledgement}
\theoremstyle{remark}
\newtheorem*{rem}{Remark}

\theoremstyle{definition}
\newtheorem{ex}{Example}[section]

\numberwithin{equation}{section}
\numberwithin{figure}{section}

\def\Hom{{\text{\rm{Hom}}}}

\def\rchi{{\hbox{\raise1.5pt\hbox{$\chi$}}}}
\def\Aut{{\text{\rm{Aut}}}}

\def\isom{\cong}

\def\dsum{\oplus}

\def\reg{{\text{\rm{reg}}}}
\def\Ghat{\hat{G}}
\def\lam{\lambda}

\def\Lieg{{\mathfrak{g}}}

\begin{document}
\large

\title[Volume of representation
varieties]
{Volume of representation
varieties}
\author[M.~Mulase]{Motohico Mulase$^1$}
\address{
Department of Mathematics\\
University of California\\
Davis, CA 95616--8633}
\email{mulase@math.ucdavis.edu}
\author[M.~Penkava]{Michael Penkava$^2$}
\address{
Department of Mathematics\\
University of Wisconsin\\
Eau Claire, WI 54702--4004}
\email{penkavmr@uwec.edu}
\date{December 1, 2002}
\thanks{$^1$Research supported   by 
NSF grant DMS-9971371
and the University of California,  Davis.}
\thanks{$^2$Research supported   by NSF grant  
DMS-0200669 and the University of Wisconsin, Eau Claire.}

\setcounter{section}{-1}
\begin{abstract}
We introduce the notion of volume 
of the representation variety of a finitely presented discrete group
in a compact Lie group  using the push-forward
measure associated to a map
defined by a presentation of the discrete group. 
We show that the
volume thus defined is invariant under the Andrews-Curtis moves
of the generators and relators of the discrete group, and moreover,  
 that  it  is actually independent of the choice of
 presentation if the difference 
of the number of 
generators and  the number of relators remains the same.
We then calculate the volume of the representation 
variety of a
surface group in an arbitrary compact Lie
group using the classical 
technique of Frobenius and Schur 
on finite groups. 
Our formulas recover
the results of Witten and Liu on the symplectic volume
and the Reidemeister torsion of the moduli space of
flat $G$-connections on a surface up to a constant factor
when the Lie group $G$ is 
semisimple. 
\end{abstract}
\maketitle

\tableofcontents
\allowdisplaybreaks

\section{Introduction and the Main Results}
\label{sect:intro}

Let $\Pi$ be a finitely presented discrete group
generated by $k$ elements with $r$ relations:
\begin{equation}
\label{eq:pi}
\Pi = \langle a_1,\dots,a_k\;|\;q_1,\dots,q_r \rangle\;.
\end{equation}
The space $\Hom(\Pi,G)$
of homomorphisms from $\Pi$ into a
compact Lie group $G$ is called the
\emph{representation variety} of $\Pi$ in $G$. 
Representation varieties of the fundamental group of 
a manifold
appear in various places in low-dimensional geometry and topology 
(see for example, \cite{CS, G,KM, Kuperberg, PX}). 
To the presentation (\ref{eq:pi})  we
associate a \emph{presentation map}
\begin{equation}
\label{eq:presentationmap}
q:G^k\owns x = (x_1,\dots,x_k)\longmapsto
(q_1(x),\dots,q_r(x))\in G^r\;,
\end{equation}
which is a real analytic map. By definition, 
we have a canonical 
identification
\begin{equation}
\label{eq:homisqinverse}
\Hom(\Pi, G) = q^{-1}(1,\dots,1)
\end{equation}
that gives a realization of the representation variety as a
real analytic subvariety of $G^k$. If the 
relators are given without redundancy, then its dimension is expected to be
\begin{equation}
\label{eq:dimofrepvar}
\dim \Hom(\Pi,G) = (k-r)\dim G\;.
\end{equation}

The purpose of this paper is to give a definition of
the \emph{volume} $|\Hom(\Pi,G)|$ of the representation variety
and to study some of its properties. 
Since we use the push-forward
measure via the presentation map (\ref{eq:presentationmap}) to
define the volume, the central question is
to determine whether it is an invariant of the group. 
In this article we give an affirmative answer to this question.
The main results of the paper are the following.

\begin{thm}
The volume is invariant under the
Andrews-Curtis moves.
\end{thm}

\noindent
Because of the dimension formula of (\ref{eq:dimofrepvar}), 
it is necessary that the difference $k-r$ 
of the numbers of generators and relators remains invariant
for different presentations of the group $\Pi$ to give the same
volume $|\Hom(\Pi,G)|$.
 Surprisingly, the condition is sufficient.

\begin{thm}
The volume is independent of the choice of  presentation
if the  difference between
the number of generators
and the number of relators is the same. 
\end{thm}

\noindent
The volume of the representation variety
of the fundamental group of a closed surface
is a topological invariant and actually expressible in terms
of irreducible representations of the compact Lie group $G$.

\begin{thm}
\label{thm:introori}
Let $S$ be a closed orientable
surface of genus $g\ge 2$. Then the volume of 
the representation variety can be computed by the 
formula
\begin{equation}
\label{eq:introvolumeori}
\sum_{\lam\in\Ghat} (\dim\lam)^{\rchi(S)}
= |G|^{\rchi(S)-1}|\Hom(\pi_1(S),G)|
\end{equation}
if the sum of the LHS converges, 
where $\Ghat$ denotes the set of isomorphism classes  of 
complex irreducible representations
of $G$, $\dim\lam$ is the complex dimension of
the irreducible representation $\lam$, and 
$\rchi(S)=2-2g$ is the Euler characteristic of $S$. 
Since $G$ is compact, $\Ghat$ is a countable set.
\end{thm}

\noindent
For a non-orientable surface $S$, the formula for the volume
involves more detailed information on the
 irreducible representations
of $G$. 
Using  the Frobenius-Schur
indicator of irreducible characters \cite{FS}, 
we decompose the set of complex irreducible representations $\Ghat$
into the union of three disjoint subsets, corresponding to
real, complex, and quaternionic irreducible representations:
\begin{equation}
\label{eq:FSindicator}
\begin{split}
\Ghat _1 &= \bigg\{\lam\in \Ghat \;\bigg|\; \frac{1}{|G|}\int_{G}
\rchi_\lam(w^2) dw = 1\bigg\}\; ;\\
\Ghat _2 &= \bigg\{\lam\in \Ghat \;\bigg|\; \frac{1}{|G|}\int_{G}
\rchi_\lam(w^2) dw = 0\bigg\}\; ;\\
\Ghat _4 &= \bigg\{\lam\in \Ghat \;\bigg|\; \frac{1}{|G|}\int_{G}
\rchi_\lam(w^2) dw = -1\bigg\}\; ,
\end{split}
\end{equation}
where $\rchi_\lam$ is the irreducible character of 
$\lam\in\Ghat$. Now we have

\begin{thm}
\label{thm:intrononori}
Let  $S$ a closed non-orientable surface.
Then we have
\begin{equation}
\label{eq:introvolumenonori}
\sum_{\lam\in\Ghat_1} (\dim\lam)^{\rchi(S)}
+ \sum_{\lam\in\Ghat_4} (-\dim\lam)^{\rchi(S)}
=  |G|^{\rchi(S)-1} |\Hom(\pi_1(S), G)|
\end{equation}
if the sum of the LHS  is absolutely convergent.
\end{thm}

The formulas (\ref{eq:introvolumeori}) and
(\ref{eq:introvolumenonori}) for a finite group,
in which case the closed surface can be anything, are
due to  Mednykh
\cite{Med} and Frobenius-Schur \cite{FS},
respectively. We refer to the excellent review article
\cite{Jones} for the history of  related topics.
For a finite group $G$ and $S=S^2$, (\ref{eq:introvolumeori}) 
reduces to the classical formula for the order of the group
\begin{equation}
\label{eq:orderformula}
|G| = \sum_{\lam\in\Ghat} (\dim \lam)^2 \;.
\end{equation}
 For a compact connected semisimple Lie group $G$, 
the formulas agree with the results of 
Witten \cite{W1991, W1992} and Liu \cite{L1,L2,L3}
up to a constant factor.
 They calculated  the natural symplectic volume 
and the Reidemeister torsion of the moduli 
space of flat $G$-connections on a surface $S$
using techniques from theoretical physics and 
analysis of heat kernels. 
The idea of using Chern-Simons gauge
theory to obtain (\ref{eq:introvolumeori})
 for a finite group was carried out
in \cite{FQ}.
The representation variety has algebro-geometric interest
as the moduli space of holomorphic vector bundles on a
compact Riemann surface, and through this interpretation
it has deep connections with conformal field theory
and topological quantum field theory. 
Let $S$ be a compact Riemann surface with a complex structure. 
Narasimhan-Seshadri theory \cite{NS1, NS2} 
shows that the moduli space
of stable holomorphic $G_\mathbb{C}$ bundles over $S$
as a real analytic variety
can be identified with
$$
\mathcal{M}(S,G_\mathbb{C})=
\frac{\Hom(\pi_1(S),G)}{G/Z(G)}\;,
$$
where $Z(G)$ is the center of $G$ and the $G$-action on 
$\Hom(\pi_1(S),G)$ is through conjugation.  
The first Chern class of the
determinant line bundle over the moduli space gives
an integral symplectic form $\omega$ on the moduli space. 
Using this $2$-form Witten \cite{W1991} computed  the natural
symplectic volume of $\mathcal{M}(S,G_\mathbb{C})$. 
His formula and ours differ by a factor of
$(2\pi)^{\dim_\mathbb{R} \mathcal{M}}$, which is due to the normalization
of the integral cohomology class $[\omega]$ on the moduli space,
but the factor does not come from the choice of the Haar measure on
the group $G$. 
From the algebro-geometric point of view, certainly the quotient of the
natural conjugate action of the representation variety is a
desirable object of study. However, our formulas
(\ref{eq:introvolumeori}) and (\ref{eq:introvolumenonori})
suggest that the volume of the
whole representation variety measured in comparison with
the volume of the space $G^{1-\rchi(S)}$ may be a more natural object.
The heat kernel method is necessary to connect the volume
of the representation variety of a surface group 
we discuss in this paper to the symplectic volume
or the Reidemeister torsion of the moduli space of 
flat connections on a closed surface. Our point is that if we
consider the natural
quantity $|G|^{\rchi(S)-1}\big|\Hom\big(\pi_1(S),G\big)\big|$,
which is independent of the
choice of the measure on $G$,  then we
can compute its value by exactly the same classical method 
of Frobenius \cite{F} and Frobenius-Schur \cite{FS} for finite groups.

The paper is organized as follows. In Section~\ref{sect:presentation}
we recall that a group $\Pi$ is completely determined by
the hom functor $\Hom(\Pi,\bullet)$. This implies that 
two presentations define the same group if and only if what we call
the presentation functors are naturally isomorphic.
With this preparation we define the volume of the
representation variety $\Hom(\Pi,G)$ in Section~\ref{sect:volume}
using a presentation of $\Pi$. We then prove that the volume
is actually independent of the presentation, utilizing the
behavior of the volume under the Andrews-Curtis moves. 
Then in Sections~\ref{sect:ori} and \ref{sect:nonori}
we calculate the volume of the representation variety of 
a surface group, the orientable case first followed by the non-orientable
case.  Throughout the paper except for Section~\ref{sect:presentation}
we need non-commutative harmonic analysis of 
distributions on a compact Lie
group and character theory. These techniques are summarized 
in Section~\ref{sect:harmonic}. The final section is devoted to
presenting a heuristic argument that relates the computation 
of the volume of this
paper to the non-commutative matrix integral 
studied in  \cite{MY,MY2}.

\begin{ack}
The authors thank Bill Goldman for drawing their
attention to important references on the 
subject of this article. They also thank 
Greg Kuperberg, Abby Thompson and Bill 
Thurston for stimulating discussions on representation
varieties and Andrews-Curtis moves.
\end{ack}

\bigskip

\section{Presentation of a Group}
\label{sect:presentation}

Consider a set of \emph{alphabets} $\{a_1,\dots,a_k\}$
consisting of $k$ letters. 
A \emph{word} $W(a)$ of length
$\nu$ is the expression of the form
\begin{equation}
\label{eq:wordform}
W(a) = a_{w(1)} ^{\epsilon(1)} \cdots a_{w(\nu)} ^{\epsilon(\nu)}\;,
\end{equation}
where the subscript
\begin{equation}
\label{eq:ellthletter}
w: \{1,\dots,\nu\}\longrightarrow \{1,\dots,k\}
\end{equation}
indicates the $\ell$-th letter $a_{w(\ell)}$ in the word $W(a)$, and
the exponent  $\epsilon(\ell)=\pm 1$ represents either the letter 
$a_{w(\ell)}$ or its \emph{inverse} $a_{w(\ell)} ^{-1}$.  
Instead of the commonly used  expression  $a_i ^\nu$, we always use
$$
\overset{\nu}{\overbrace{a_i \cdots a_i}}\;,
$$
which is a word of length $\nu$. 
Throughout this article,  every word is assumed to be reduced,
namely,  no expressions like $a_i a_i ^{-1}$ and $a_i ^{-1} a_i$
appear in a word. 

For every given group $G$, 
a word $W(a)$ of letters $a_1,\dots,a_k$ (that are unrelated
with $G$) defines a map
\begin{equation}
\label{eq:wordmap}
W:G^k \owns x=(x_1,\dots, x_k)\longmapsto
W(x) = x_{w(1)} ^{\epsilon(1)} \cdots x_{w(\nu)} ^{\epsilon(\nu)}
\in G\;.
\end{equation}
 The map $W$ is the composition of two maps
\begin{equation}
\label{eq:wordmapfactor}
\begin{CD}
G^k @. \quad \owns\quad  @. (x_1,\dots,x_k)\\
@V{\iota_W}VV @. @VV{\iota_W}V\\
G^\nu @. \quad \owns\quad @.
 (x_{w(1)} ^{\epsilon(1)}, \dots ,x_{w(\nu)} ^{\epsilon(\nu)})\\
@V{m}VV @. @VV{m}V\\
G @. \quad \owns\quad @. 
x_{w(1)} ^{\epsilon(1)} \cdots x_{w(\nu)} ^{\epsilon(\nu)}\;,
\end{CD}
\end{equation}
where $\iota_W:G^k\rightarrow G^\nu$ represents the
shape of the word $W(a)$,
and 
\begin{equation}
\label{eq:multiplicationmap}
m:G^\nu \owns (y_1,\dots,y_\nu)\longmapsto 
y_1\cdots y_\nu\in  G
\end{equation}
is the multiplication map. If the word $W(a)$ does not contain
the letter $a_i$, then $x_i$ is sent to $1\in G$ by the map $\iota_W$.
In particular, if $W=\emptyset$ is the empty word, then
$\iota_\emptyset$ is the trivial map that
sends  everything in $G^k$ to $\{1\}
= G^0$.

Let us consider a finitely presented discrete group
\begin{equation}
\label{eq:grouppres}
\Pi = \langle a_1,\dots,a_k\;|\;
q_1,\dots,q_r\rangle
\end{equation}
defined  by $k$ generators and $r$ relations $q_1=\cdots =q_r = 1$. 
Every \emph{relator} $q_j$ is a reduced word, and we always assume that
none of the relators is the empty word. 
The set of relators
defines the \emph{presentation map}
\begin{equation}
\label{eq:mapq}
q:G^k\owns x= (x_1,\dots,x_k)\longmapsto
(q_1(x),\dots,q_r(x))\in G^r\;.
\end{equation}
We give a group structure to  the target space $G^r$ as the product group
of $r$ copies of $G$. The identity element of $G^r$ is
$1 = (1,\dots,1)$. 
 The \emph{representation space} of $\Pi$
in $G$ is the set of group homomorphisms from $\Pi$ to $G$
and is naturally identified with
\begin{equation}
\label{eq:repvar}
\Hom(\Pi, G) = q^{-1}(1)\;.
\end{equation}

Let $\mathcal{G}$ denote the category of groups with
group homomorphisms as  morphisms among the objects. A
finite presentation (\ref{eq:grouppres}) defines a 
covariant functor from
$\mathcal{G}$ to the category $\mathcal{S}$ of sets,
associating $q^{-1}(1)$ to a group $G$. Let us call this 
functor the \emph{presentation functor}.
 The covariance of the presentation functor is a consequence of the 
commutativity of the diagram 
\begin{equation}
\label{eq:covariance}
\begin{CD}
H^k @>q>> H^r\\
@V{\gamma}VV @VV{\gamma}V\\
G^k@>>q> G^r
\end{CD}
\end{equation}
for every homomorphism 
$\gamma:H\rightarrow G$,
where the vertical arrows are the natural maps induced by
 $\gamma$. This functor is representable and
the group $\Pi$ is the representing object in $\mathcal{G}$
\cite{Mac}. 
Thus the presentation functor, whose representation 
is the hom functor
$\Hom(\Pi,\bullet)$, is determined by the group $\Pi$ alone and is
independent of the choice of its presentation. 
Let $\mu:\Pi\rightarrow \Pi'$ be a group homomorphism
from $\Pi$ to $\Pi'$. The pull-back map $\mu^*$ defines
the natural morphism
$$
\mu^*:\Hom(\Pi',\bullet)\longrightarrow \Hom(\Pi,\bullet)
$$
between the hom functors. 

Notice that the functor
from $\mathcal{G}$ to the functor category 
$\mathcal{F}(\mathcal{G},\mathcal{S})$ that associates
$$
\mathcal{G}\owns \Pi \longmapsto 
\Hom(\Pi,\bullet)\in \mathcal{F}(\mathcal{G},\mathcal{S})
$$
is a faithful contravariant
functor. Indeed, if  there is a natural isomorphism
$$
A: \Hom(\Pi,\bullet)\overset{\sim}{\longrightarrow} 
\Hom(\Pi',\bullet)\;, 
$$
then by applying this
functor to $\Pi'\in \mathcal{G}$, we obtain 
 a homomorphism $\mu:\Pi\rightarrow \Pi'$
that corresponds  to the identity automorphism  
$1_{\Pi'}\in \Hom(\Pi',\Pi')$. Similarly, $1_{\Pi}\in \Hom(\Pi,\Pi)$
corresponds to $\xi:\Pi'\rightarrow \Pi$. Because of the naturalness 
condition, for the homomorphism $\mu$ we have the commutative
diagram
\begin{equation*}
\begin{CD}
\Pi @>{1_\Pi}>> \Pi @.  \hskip0.5in
\Hom(\Pi, \Pi) @>{A_\Pi}>>\Hom(\Pi',\Pi)
@. \hskip0.5in \Pi @<{\xi}<< \Pi'\\
@.  @V{\mu}VV \hskip0.5in
 @VVV @VVV
\hskip0.5in @V{\mu}VV @.\\
@. \Pi' @.\hskip0.5in
\Hom(\Pi,\Pi') @>{A_{\Pi'}}>> \Hom(\Pi',\Pi')
@. \hskip0.5in \Pi'
\end{CD}
\end{equation*}
that relates 
\begin{equation*}
\begin{CD}
 1_{\Pi} @>>> \xi\\
@VVV @VVV\\
\mu @>>> 1_{\Pi'}
\end{CD}\;.
\end{equation*}
In particular, we have $\mu\circ \xi = 1_{\Pi'}$, hence $\Pi\isom \Pi'$. 
As a consequence of this faithfulness,
we note  that  if two presentation functors are
naturally isomorphic, then  the two presented groups
are actually isomorphic.

The above consideration suggests that we can 
use a particular presentation of a group $\Pi$ to define 
the  \emph{volume} 
of the representation space $\Hom(\Pi,G)$, 
which will be an invariant of 
the group $\Pi$, independent of its presentation. 
However,
to carry out the calculation of
 the volume  $|\Hom(\Pi,G)|$, we need not only 
 a presentation of    $\Pi$, but also
 a more specific natural isomorphism between
presentation functors. 
Notice that two presentations 
\begin{equation}
\label{eq:twopres}
\Pi_1=\langle a_1,\dots,a_k\;|\;q_1\dots,q_r\rangle
\qquad\text{and} \qquad 
\Pi_2=\langle b_1,\dots,b_k\;|\;s_1\dots,s_r\rangle
\end{equation}
define the same group
 if and only if for every $a_i$ there is a word
$a_i(b)$ of the generators $b_1,\dots,b_k$ and for
every $b_j$ there is a word $b_j(a)$ of $a_1,\dots,a_k$
such that
for any
group $G\in \mathcal{G}$, 
the maps $a$ and $b$ associated to these words are bijective and
satisfy the commutativity 
\begin{equation}
\label{eq:groupequivalenceCD}
\begin{CD}
G^k @= G^k @>q>> G^r\\
@A{a}A{\wr}A @V{\wr}V{b}V\\
G^k @= G^k @>>s> G^r
\end{CD}
\end{equation}
and
\begin{equation}
\label{eq:groupequivalence}
q^{-1}(1) = b^{-1}\big(s^{-1}(1)\big),
\qquad 
\qquad 
s^{-1}(1) = a^{-1}\big(q^{-1}(1)\big)\;.
\end{equation}
 In Section~\ref{sect:volume} we use this specific form
of the natural isomorphism of the presentation functors
to show the independence of the volume of the
representation variety on the choice of presentation.

\bigskip

\section{Volume of Representation Varieties}
\label{sect:volume}

From now on we restrict our attention to the category of 
compact real analytic Lie groups. The 
representation space $\Hom(\Pi, G)$ of 
$$
\Pi = \langle a_1,\dots,a_k\;|\;q_1,\dots,q_r\rangle
$$
 in a 
compact Lie group $G$ is called the \emph{representation variety}. 
Unless otherwise stated, we do not assume connectivity or
semisimplicity for the compact Lie group $G$. Therefore, it
can be in particular a finite group. The natural identification 
(\ref{eq:repvar}) makes $\Hom(\Pi, G)$ a real analytic
subvariety of $G^r$. In this section we define the volume of
the representation variety using techniques from non-commutative
harmonic analysis on compact Lie groups. 
Especially we use theory of 
distributions on a compact Lie group and character theory 
freely in Sections~\ref{sect:volume}, \ref{sect:ori} and
\ref{sect:nonori}. We refer to the summary presented 
in Section~\ref{sect:harmonic} for terminologies,  definitions and
properties of distributions necessary for our investigation.

We denote by $dx$ a left and right invariant measure on $G$. 
We do not normalize the invariant measure, and call the quantity
$$
|G|=\int_G dx
$$
the volume of the group $G$. 
The invariant measure on the product group $G^r$ is the
product measure $dw_1\cdots dw_r$.
The $\delta$-function on $G^r$
is defined by 
\begin{equation}
\label{eq:multidelta}
\delta_r(w_1,\dots,w_r) = \delta(w_1)\cdots\delta(w_r)\;.
\end{equation}
For the presentation map
$$
q:G^k\longrightarrow G^r
$$
of (\ref{eq:mapq}), we associate the
\emph{volume distribution} $f_q(w)$ on $G^r$  
by
\begin{equation}
\label{eq:volumedist}
f_q(w) = |q^{-1}(w)| = 
\int_{G^k} \delta_r(q\cdot w^{-1})dx_1\cdots dx_k\;,
\end{equation}
where 
$$
q\cdot w^{-1} = (q_1(x) w_1 ^{-1},\dots,q_r(x)w_r ^{-1})
\in G^r
$$
and $dx_1\cdots dx_k$ is the product measure on $G^k$.
The volume distribution is a Schwartz distribution on the
group $G^r$ characterized by the continuous linear functional
associated with the push-forward measure
\begin{equation}
\label{eq:fqcharacterization}
g(w)\longmapsto
\int_{G^r} f_q(w) g(w)dw_1\cdots dw_r
 = \int_{G^k} g(q(x_1,\dots,x_k))
dx_1\cdots dx_k
\end{equation}
for every $g(w)\in C^\infty(G^r)$.
Sard's theorem tells us that the set of critical values of
$q$ has measure $0$ in the target space  with respect to the
Lebesgue measure of $G$ as a real analytic manifold, and hence also 
 with respect to the Haar measure $dx$. 
If the distribution
$f_q(w)$ is regular at $w=1$ 
(see Section~\ref{sect:harmonic}), then we define
the \emph{volume} of the representation variety by
\begin{equation}
\label{eq:defofvolume}
|\Hom(\Pi,G)| = f_q(1)\;.
\end{equation}
The total volume of $G^k=q^{-1}(G^r)$
with respect to the volume distribution is
$$
\int_{G^r} f_q(w) dw_1\cdots dw_r = 
\int_{G^k}\left(\int_{G^r} \delta_r(q\cdot w^{-1}) dw_1\cdots dw_r
\right)dx_1\cdots dx_k = |G^k|\;.
$$
Thus our definition (\ref{eq:defofvolume}) agrees with Fubini's volume
if the presentation map $q$ is a fiber bundle.
Since evaluation of a distribution at a point does not usually
make sense, we caution that the volume is defined \emph{only when}
$f_q(1)$ is regular. This imposes a strong condition on the
presentation of $\Pi$.  

\begin{prop}
\label{prop:fqregularity}
If the presentation map $q:G^k\rightarrow G^r$ has $1\in G^r$ as a
regular value, then the volume distribution $f_q(w)$ is also regular
at the identity element.
\end{prop}

\begin{proof}
The condition implies that the differential $dq$
of the presentation map $q$ has the maximal rank at every point of 
$q^{-1}(1)$, and that $q^{-1}(1)$ is a non-singular
submanifold of $G^k$.
 Since $G$ is compact, there is an open neighborhood
$U$ of $1\in G^r$ such that 
\begin{equation}
\label{eq:phi}
\phi: F\times U \overset{\sim}{\longrightarrow}q^{-1}(U)   
\end{equation}
as a $C^\infty$ manifold. Here $F=\Hom(\Pi,G)=q^{-1}(1)$
 is a $C^\infty$ manifold
and the presentation map restricted to $q^{-1}(U)$ is the projection 
$F\times U \rightarrow U$. Let $\omega$ be the invariant volume
form on $G$ corresponding to the invariant measure $dx$. Denote by
$p_j:G^k\rightarrow G$ the projection to the $j$-th component.
Then 
\begin{equation}
\label{eq:omegak}
\omega^k \overset{\text{def}}{=}
 p_1 ^*\omega \wedge \cdots \wedge p_k ^* \omega 
\end{equation}
corresponds to the 
product measure $dx_1\cdots dx_k$ on $G^k$. 
Similarly, we denote by $\omega^r$ the invariant volume form on 
the product group $G^r$ corresponding to $dw_1\cdots dw_r$.
 Notice that the volume form $\phi^*\omega^k$ on 
$F\times U$ is a $C^\infty$ form. Therefore, there is a $C^\infty$
function $h(w)$ on $U$ such that the integration of
$\phi^*\omega^k$ along fiber associated with the fibration
$F\times U \rightarrow U$ is
given by
$$
\int_{F\times \{w\}} \phi^*\omega^k = h(w)\; \omega^r\;. 
$$
Since the volume distribution $f_q(w)$ is the push-forward measure, 
we have the equality
$$
\big(f_q|_U\big)(w) = h(w)
$$
as a distribution. Therefore, 
the localization  $f_q|_U$ of the volume distribution $f_q(w)$ to $U$ is
a smooth function.
\end{proof}

\begin{ex}
Consider a presentation
$$
\Pi = \langle a,b\;|\;ab^{-1}=1\rangle\;.
$$
The corresponding presentation map is
$$
q: G\times G\owns (x,y)\longmapsto xy^{-1}\in G\;,
$$
hence $q^{-1}(1)$ is the diagonal of $G\times G$. 
The volume of the diagonal is calculated by
$$
|\Hom(\Pi, G)| =\int_{G^2}\delta(xy^{-1})dxdy = |G|\;.
$$
Notice that the same group is presented by $\langle a\rangle$. 
With this presentation, there is no relation, and the presentation map
is the trivial map
$$
G\longrightarrow \{1\}=G^0\;.
$$
The delta function on the trivial group is the constant $1\in\mathbb{C}$,
hence 
$$
|\Hom(\langle a\rangle,G)|=\int_G dx = |G|\;.
$$
\end{ex}

For a compact Lie group $G$, we require that
the maps $a$ and $b$ of the natural isomorphism of
(\ref{eq:groupequivalenceCD}) and (\ref{eq:groupequivalence})
are real analytic automorphisms of the real
analytic manifold $G^k$.

\begin{lem}
\label{lem:volumeforminv}
Suppose 
$$
b:G^k\owns x= (x_1,\dots,x_k)\longmapsto 
(b_1(x),\dots,b_k(x)) = b(x) \in G^k
$$ 
is an analytic automorphism of the real analytic maniforld $G^k$
given by $k$ words $b_1,\dots,b_k$ of $x_1,\dots,x_k \in G$
 such that its inverse  is
also given by $k$ words in a similar way. Then the volume
form $\omega^k$ of $G^k$ is invariant under the automorphism
$b$ up to a sign:
\begin{equation}
\label{eq:binvarianceofomegak}
b^*\omega^k = \pm \omega^k\;.
\end{equation}
\end{lem}

\begin{proof}
Since the automorphism $b$ maps $(1\dots,1)\in G^k$
to itself, the differential $db$ at $(1\dots,1)\in G^k$ can be described
in terms of  the Lie algebra $\Lieg$ of $G$. 
Noticing that a word is a product of $\pm 1$ powers of $x_j$'s, the 
differential 
$$
db:T_{(1,\dots,1)} G^k = \Lieg ^{\dsum k}
\overset{\sim}{\longrightarrow}
 \Lieg ^{\dsum k} = T_{(1,\dots,1)} G^k
$$
is given by an invertible $k\times k$ matrix of integer entries $N\in 
GL(k,\mathbb{Z})$.  Since $\det N = \pm 1$, the $G$-invariant
volume form on $\mathfrak{g}^{\dsum k}$ is preserved by
$db$ up to the sign $(\det N)^{k\dim G}$. 

Now notice that 
\begin{equation}
\label{eq:omegakcomparizon}
\begin{split}
b^* \omega^k &= b_1 ^*\omega \wedge \cdots \wedge b_k ^*\omega\\
\omega^k &= p_1 ^*\omega \wedge \cdots \wedge p_k ^*\omega\;,
\end{split}
\end{equation}
where $p_j:G^k\rightarrow G$ is the projection to the $j$-th
component used in (\ref{eq:omegak}). 
 To analyze the pull-back of $\omega$ by a word
$b_i$, let us examine the pull-back of $\omega$
via the multiplication map
$$
m:G^\nu \owns (y_1,\dots,y_\nu)\longmapsto
y_1\cdots y_\nu\in  G
$$ 
of (\ref{eq:multiplicationmap}). Here $\nu$ is the length
of  $b_j$.
Since $\omega$ is left and right multiplication invariant, 
the restriction of $m^* \omega$ to the $\ell$-th factor of $G^\nu$ is
$\omega$ itself. To be more precise, it is $p_\ell ^* \omega$. 
If the $\ell$-th letter of the word $b_i$ is $x_1$, then the 
restriction of 
$$
\iota_{b_i}:G^k \longrightarrow G^\nu
$$ 
of (\ref{eq:wordmapfactor})
on the first factor of $G^k$ into the $\ell$-th
factor of $G^\nu$ is the identity map. Thus 
$\iota_{b_i} ^* m^* \omega$ restricted to the first factor
of $G^k$ is just an integer multiple of $\omega$ (or to be
precise $p_1 ^* \omega$), and the coefficient is the number
of times the letter $x_1$ appears in the word $b_i$ plus
$(-1)^{\dim G}$ for each appearance of  $x_1 ^{-1}$
in $b_i$. Define $n_{ij}$ as the number of letter $x_j$
in the word $b_i$ minus the number of $x_j ^{-1}$ in $b_i$. 
Then the matrix $N = (n_{ij})$ is the same matrix $N$ that
represents the differential $db$. Therefore, 
$N\in GL(k,\mathbb{Z})$. Noticing the identity
$$
 b^* \omega ^k = 
b_1 ^*\omega \wedge \cdots \wedge b_k ^*\omega
=\iota_{b_1} ^* m^* \omega\wedge \cdots \wedge
\iota_{b_k} ^* m^* \omega\;,
$$
we obtain the desired formula
$$
b^* \omega ^k = N^* \omega^k  = (\det N)^{k\dim G}\cdot \omega ^k\;.
$$
\end{proof}

A consequence of this lemma is the independence of
the volume $|\Hom(\Pi,G)|$ on the choice of presentation
among those with $k$ generators and 
$r$ relators.

\begin{thm}
\label{thm:indeponpresentation}
Let $G$ be a compact Lie group.
Suppose we have two presentations (\ref{eq:twopres})
of the same group $\Pi$ with the same $k$ and $r$.
If the presentation maps $p$ and $q$  have the regular value
at $1\in G^r$, then the volume of the representation variety of 
$\Pi$ in $G$ is independent of the presentation:
\begin{equation}
\label{eq:fq=fr}
f_q(1)=\int_{G^k} \delta_r\big(q(x_1,\dots,x_r)\big)dx_1\cdots dx_k
= \int_{G^k} \delta_r\big(s(y_1,\dots,y_k)\big)dy_1\cdots dy_k
=f_s(1)\;.
\end{equation}
\end{thm}

\begin{proof}
As in Proposition~\ref{prop:fqregularity}, there is an open neighborhood
$U\owns 1$ of $G^r$ and $C^\infty$ maps $\phi$ and $\psi$
such that
\begin{equation*}
\begin{CD}
q^{-1}(U) @<\sim<\phi< F\times U @>\sim>\psi> (sb)^{-1}(U)\\
@VqVV @VV{\text{proj}}V @VVsbV\\
U @= U @= U
\end{CD}\;,
\end{equation*}
where $F=q^{-1}(1) = (sb)^{-1}(1)\subset G^k$. Note in particular that
the restrictions of $\phi$ and $\psi$  on $F\times \{1\}$
are the identity map.
Let $\omega^k$ and $\omega^r$ denote the volume forms
on $G^k$ and $G^r$, respectively, defined by
(\ref{eq:omegak}). Then the localizations  of the
volume distributions $f_q$ and 
$f_{sb}$ on $U$ (see Section~\ref{sect:harmonic})
are  $C^\infty$ functions
and given by the integration along fiber:
$$
\int_{F\times \{w\}} \phi^*\omega^k  = \big(f_q|_U\big)(w)\;\omega^r\;,
\qquad\qquad 
\int_{F\times \{w\}} \psi^*\omega^k = \big(f_{sb}|_U\big)(w)\;\omega^r\;.
$$
Choose a volume form $\omega_F$ on $F$ and write
$$
\phi^*\omega^k = h_q(z,w) \omega_F\wedge \omega^r
\qquad \text{and} \qquad
\psi^*\omega^k = h_{sb}(z,w) \omega_F\wedge \omega^r
$$
using $C^\infty$ functions $h_q$ and $h_{sb}$ on $F\times U$. 
Then the integration along fiber is given by
\begin{equation*}
\begin{split}
\big(f_q|_U\big)(w)\;\omega^r&=\int_{F\times \{w\}} \phi^*\omega^k
=\left(\int_{F\times \{w\}} h_q(z,w) \omega_F\right) \; \omega^r\\
\big(f_{sb}|_U\big)(w)\;\omega^r&=\int_{F\times \{w\}} \psi^*\omega^k
=\left(\int_{F\times \{w\}} h_{sb}(z,w) \omega_F\right) \;\omega^r\;.
\end{split}
\end{equation*}
Since $h_q(z,1)=h_{sb}(z,1)$ as a function on $F\times \{1\}$, 
we obtain $f_q(1) = f_{sb}(1)$. 
Thus if we
have $b^*\omega^k = \omega^k$, then 
$f_{sb}(1) = f_s(1)$, and the proof is completed. 
This invariance of the volume form $\omega^k$ 
under the automorphism $b$ is assured by the previous lemma.
The sign $\pm 1$ in Lemma~\ref{lem:volumeforminv}
is not our concern here because we
 always use the positive volume form for integration on $G^k$.
\end{proof}

To establish the independence of the volume on 
the choice of presentation with different numbers of 
generators and relators, we first need its  
invariance under the {Andrews-Curtis moves}
\cite{AC}.

\begin{Def}
\label{Def:AC}
The \emph{Andrews-Curtis moves} of a presentation 
(\ref{eq:grouppres}) of
a discrete group  are the following six  
operations on the generators and relators:
\begin{enumerate}
\item Interchange $q_1$ and $q_j$ for $j=2,\dots,r$.
\item Replace $q_1$ with $a q_1 a^{-1}$ for an element
$a\in \Pi$.
\item Replace $q_1$ with $q_1 ^{-1}$.
\item Replace $q_1$ with $q_1 q_2$.
\item Add another generator $a$ and a relator
$q_{r+1}= a$. 
\item Delete a generator $a$ and a relator $a$.
\end{enumerate}
\end{Def}

\begin{thm}
\label{prop:ACinvariance}
Let $\Pi = \langle a_1,\dots, a_k\;|\; q_1,\dots, q_r\rangle$ be
a finite presentation of a discrete group and $G$ a compact Lie
group. If the volume distribution $f_q(w)$ is
regular at $w=1\in G^r$ and hence the volume 
of the representation variety $f_q(1)
= |\Hom(\Pi, G)|$ is well-defined, then the volume is 
invariant under the  Andrews-Curtis moves. 
\end{thm}

\begin{proof}
The invariance under the first three moves follows immediately 
from (\ref{eq:multidelta}) and (\ref{eq:formulasofdelta}). 
To show the invariance under (4), we note the following identity 
of $\delta$-functions:
\begin{equation}
\label{eq:2deltaidentity}
\delta_2(x,y) = \delta(x)\delta(y) = \delta(xy)\delta(y)\;.
\end{equation} 
This is because for every $g(x,y)\in C^\infty(G^2)$, we have
\begin{equation*}
\begin{split}
\int_{G^2} g(x,y)\delta(xy)\delta(y)dxdy
&=\int_G g(y^{-1},y)\delta(y)dy\\
&=g(1,1) \\
&= \int_{G^2} g(x,y)\delta_2(x,y)dxdy\;.
\end{split}
\end{equation*}
Therefore, 
$$
\delta(q_1)\delta(q_2)\cdots \delta(q_r)
=\delta(q_1q_2)\delta(q_2)\cdots \delta(q_r)\;.
$$
The invariance under (5) and (6) comes from the following
commutative diagram
\begin{equation}
\label{eq:AC56}
\begin{CD}
G^k\times G @>(q,id)>> G^r\times G\\
@A{i_k}AA @AA{i_r}A\\
G^k @>>q> G^r
\end{CD}\;,
\end{equation}
where the vertical maps are the inclusions into the 
first $G^k$ and $G^r$ components, respectively,
with $1\in G$ in the last component.
The invariance is a consequence of the identity
$$
\int_{G^k\times G} \delta_r (q)\delta(x)dx_1\cdots dx_k dx
=
\int_{G^k} \delta_r (q) dx_1\cdots dx_k\;.
$$
\end{proof}

The Andrews-Curtis moves were introduced in \cite{AC}
with the well-known conjecture stating that
 every presentation of the trivial group with the same number
of generators and relations can be transformed into the 
standard presentation by a finite sequence of moves 
(1)--(6). The conjecture is still open to date. 
For three-manifold groups,
the conjecture is expected to be true.

A three-manifold group $\Pi = \pi_1(M)$ of a closed oriented
$3$-manifold $M$ always has a \emph{balanced} presentation,
meaning that the number of generators and relations are the same,
given by a Heegaard splitting of $M$.
If $\Pi$ has a balanced presentation, then the expected
dimension of the representation variety $\Hom(\Pi,G)$ in a
compact Lie group $G$ is $0$, and hence it consists of a finite set.
The Casson invariant counts the number of elements of
$\Hom\big(\pi_1(M),SU(2)\big)$ with sign. 
If indeed $\Hom(\pi_1(M),G)$ is a finite set, then its cardinality 
$|\Hom(\pi_1(M),G)|$ is a
topological invariant of $M$. For example, Kuperberg \cite{Kuperberg}
discussed this quantity for a finite group $G$. 
But often the representation 
variety contains higher dimensional exceptional components.
If such a situation happens, then the volume distribution becomes
singular at the identity element of $G^r$, and our theory
does not provide any useful information.

We are now ready to prove that the volume of
the representation variety is
 independent of the choice
of presentation, if it is well-defined. More precisely, we have

\begin{thm}
\label{thm:indeponpresentation2}
Suppose a group $\Pi$ has two presentations of the form
\begin{equation}
\label{eq:presentationone}
\langle a_1,\dots,a_k \;|\; q_1,\dots, q_r\rangle
\end{equation}
and
\begin{equation}
\label{eq:presentationtwo}
\langle b_1,\dots, b_{\ell} \;|\; s_1,\dots, s_t\rangle
\end{equation}
such that 
\begin{equation}
\label{eq:samedifference}
k-r = \ell -t\;,
\end{equation}
or in other words, 
the expected dimensions of $\Hom(\Pi,G)$ through the two
presentations agree. If the volume distributions
$f_q$ and $f_s$ are regular at the identity element of 
the product groups $G^r$ and $G^t$, respectively, then 
$$
|\Hom(\Pi,G)| = f_q(1) = f_s(1)\;.
$$
\end{thm}

\begin{proof}
Without loss of generality, we can assume $k\le \ell$. 
Using the Andrews-Curtis moves, we supply the new
generators $a_{k+1},\dots,a_\ell$ and relators
$$
q_{r+1}=a_{k+1},\dots, q_{r+\ell -k} = a_\ell
$$
to the first presentation (\ref{eq:presentationone}) so that we have
\begin{equation}
\label{eq:enlargedpres}
\Pi = \langle a_1,\dots, a_\ell\;|\; q_1,\dots, q_{r+\ell -k}\rangle\;.
\end{equation}
From (\ref{eq:samedifference}), we know $r+\ell - k = t$. 
Therefore, the two presentations (\ref{eq:enlargedpres})
and (\ref{eq:presentationtwo}) have the same numbers of 
generators and relators. If $f_q(w)$ is regular at $w=1\in G^r$, then
the volume distribution of the enlarged presentation 
(\ref{eq:enlargedpres}) is also regular at $1\in G^\ell$. 
Therefore, from Theorem~\ref{thm:indeponpresentation}, 
we conclude that
$$
f_q(1) = f_s(1)\;.
$$
\end{proof}

For a one-relation group such as a surface group, 
character theory of compact Lie groups
provides a computational tool. In the next two sections,
we calculate the volume distribution for every surface
group in terms of irreducible characters of  $G$. 

\begin{lem}
\label{lem:class}
Let $\Pi = \langle a_1,\dots,a_k\;|\;q \rangle$
be a one-relation group. Then for every compact Lie group 
$G$, the volume distribution (\ref{eq:volumedist})
is a class distribution on $G$.
\end{lem}

\begin{proof}
For  $y\in G$ and $g(w)\in C^\infty(G)$, we have
\begin{equation*}
\begin{split}
\int_G f_q(ywy^{-1}) g(w)dw &= \int_G f_q(w) g(y^{-1}wy)dw\\
&=\int_{G^k}g(y^{-1}q(x_1,\dots,x_k) y)
dx_1\cdots dx_k\\
&=\int_{G^k}g(q(y^{-1}x_1y,\dots,y^{-1}x_k y))
dx_1\cdots dx_k\\
&=\int_{G^k}g(q(x_1,\dots,x_k))dx_1\cdots dx_k\\
&=\int_G f_q(w) g(w)dw
\end{split}
\end{equation*}
because of the left and right invariance of the Haar measure. 
\end{proof}

\noindent
Since the irreducible characters are real analytic
functions on $G$ and form an orthonormal 
basis for the $L^2$ class functions on $G$, we have
\begin{lem}
\label{lem:charexpansion}
The character expansion of the class distribution
$f_q(w)$ on $G$ is given by
\begin{equation}
\label{eq:characterexpansionoff}
f_q(w) = \frac{1}{|G|}\sum_{\lam\in\Ghat}\left(
\int_{G^k} \overline{\rchi}_\lam(q(x_1,\dots,x_k))dx_1\cdots
dx_k\right) \rchi_\lam(w)\;.
\end{equation}
\end{lem}

\noindent
Notice that the integral in (\ref{eq:characterexpansionoff}) is
convergent and well-defined, and the infinite sum is convergent
with respect to the strong topology of $\mathcal{D}'(G)$. Thus
the regularity question of $f_q(w)$ can be answered
by studying the convergence 
problem of
the infinite sum over all irreducible representations. 
For example, if the sum of (\ref{eq:characterexpansionoff}) is
convergent with respect to the sup norm on $G$,
 then $f_q$ is a $C^\infty$ function
on $G$ (see Section~\ref{sect:harmonic}).

\bigskip

\section{Representation Varieties of Orientable Surface
Groups}
\label{sect:ori}

In this section we calculate the volume distribution
of the surface group $\pi_1(S)$ in $G$ and the volume
$|\Hom(\pi_1(S),G)|$ of the representation variety
for a closed orientable 
surface $S$ and a compact Lie group $G$.
In this and  the next section, the statements and the 
argument of proof are actually almost identical to that
of the classical papers by Frobenius and Schur 
\cite{F,FS}. Because of the way we define  
the volume in the previous section,
 the same argument  simultaneously works for both finite
groups and compact Lie groups.

Let $S$ be a closed orientable surface of genus $g\ge 1$. Its
fundamental group is given by
\begin{equation}
\label{eq:pi1ori}
\pi_1(S) = \langle a_1,b_1,\dots,a_g,b_g\;|\;
[a_1,b_1]\cdots [a_g,b_g]\rangle\;,
\end{equation}
where $[a,b]=aba^{-1}b^{-1}$ is the group commutator
of $a$ and $b$. The presentation map is a real analytic map
\begin{equation}
\label{eq:qg}
q_{g}:G^{2g}\owns (x_1,y_1,\dots,x_g,y_g)
\longmapsto [x_1,y_1]\cdots [x_g,y_g]\in G\;.
\end{equation}
 The  
representation variety of  $\pi_1(S)$ in $G$ is a real analytic
subvariety 
$
\Hom(\pi_1(S),G) = q_{g} ^{-1}(1)
$
 of $G^{2g}$.
For every $g\ge 1$, we denote by $f_g$  the volume distribution
\begin{equation}
\label{eq:fg}
f_g(w) = f_{q_g}(w)
=\int_{G^{2g}} \delta\big([x_1,y_1]\cdots
[x_g,y_g] w^{-1}\big) dx_1dy_1
\cdots dx_{g}dy_g\;.
\end{equation}
By  Lemma~\ref{lem:class}, $f_g$ is a class distribution.
Our purpose is to calculate its character expansion. 
From the convolution property of the delta function
$$
\delta(q_{g+h}w^{-1})
=\int_G \delta(q_{h}w^{-1}u^{-1})\delta(uq_{g})du\;,
$$
we have
\begin{equation}
\label{eq:fgfh}
f_{g+h} = f_{g}*f_{h} = f_{h}*f_{g}
\end{equation}
for every $g,h\ge 1$. More precisely, consider the breakdown of 
$q_{g+h}$:
\begin{equation*}\begin{CD}
G^{2g}\times G^{2h}@.\quad
\owns @.  (x_1,y_1,\dots, x_{g+h},y_{g+h})\\
@V{(q_g\times q_h)}VV @. @VVV\\
G\times G @. \quad \owns @. ([x_1,y_1]\cdots [x_g,y_g], [x_{g+1},y_{g+1}]\cdots 
[x_{g+h},y_{g+h}] )\\
@V{\text{multiplication}}VV @. @VVV\\
G @. \quad\owns @. [x_1,y_1]\cdots [x_{g+h},y_{g+h}]\;.
\end{CD}\end{equation*}
The distribution $f_{g+h}$ is the 
push-forward of the product measure on $G^{2g+2h}$
via $q_{g+h}$,  which factors into the convolution product
 according to 
(\ref{eq:fgfh}). Here we have used  the fact that the convolution of two
distributions $f_g$ and $f_h$ corresponds to the push-forward
measure via the multiplication map $G\times G\rightarrow G$.  
In particular, $f_g$ is the $g$-th convolution power of $f_1$:
$$
f_g = \overset{g\text{-times}}{\overbrace{f_1 *\cdots * f_1}}
\;.
$$
Therefore,   it suffices to find the character expansion of $f_1$.
The following formula and its proof is essentially
due to Frobenius \cite{F} published in 1896. We  translate his
argument from a finite group to a compact Lie group
using the Dirac delta function and integration over the group.
 
\begin{prop}
\label{prop:f1}
As a distribution on $G$, we have the following 
strongly convergent character expansion formula in $\mathcal{D}'(G)$:
\begin{equation}
\label{eq:f1}
f_1(w)= \int_{G^2} 
\delta(xyx^{-1}y^{-1}w^{-1})dxdy
= \sum_{\lam\in\Ghat} \frac{|G|}{\dim\lam} \rchi_\lam(w)\;.
\end{equation}
\end{prop}

\begin{proof}
We begin with the convolution formula for the push-forward measure
$$
\delta(xyx^{-1}y^{-1}w^{-1}) 
= \int_G \delta(x(wu^{-1})^{-1})\delta(yx^{-1}y^{-1}u^{-1})du\;.
$$
We use the class distribution of (\ref{eq:eta})
and its character expansion (\ref{eq:expansionofeta}) here to obtain
$$
\eta_{x^{-1}}(u) = \int_G \delta(yx^{-1}y^{-1}u^{-1})dy
= \sum_{\lam\in\Ghat} \overline{\rchi}_\lam(x^{-1})
\rchi_\lam(u)
= \sum_{\lam\in\Ghat} {\rchi}_\lam(x)
\rchi_\lam(u)\;.
$$
Then we have
\begin{equation*}
\begin{split}
f_1(w)&= \int_{G^2} 
\delta(xyx^{-1}y^{-1}w^{-1})dxdy\\
&=\int_{G^3}
\delta(x(wu^{-1})^{-1})\delta(yx^{-1}y^{-1}u^{-1})dudxdy\\
&=\int_{G^2}
\delta(x(wu^{-1})^{-1})\eta_{x^{-1}}(u)dudx\\
&=\sum_{\lam\in\Ghat}\int_{G^2}
\delta(x(wu^{-1})^{-1}){\rchi}_\lam(x)
\rchi_\lam(u)dudx\\
&=\sum_{\lam\in\Ghat}\int_{G}
\rchi_\lam(wu^{-1})\rchi_\lam(u)du\\
&=\sum_{\lam\in\Ghat}\frac{|G|}{\dim\lam}\rchi_\lam(w)\;.
\end{split}
\end{equation*}
\end{proof}

Taking the $g$-th convolution power of $f_1$, we obtain

\begin{thm}
\label{thm:fg}
For every $g\ge 1$, the character expansion of the
volume distribution $f_g = f_{q_g}$ is given by
\begin{equation}
\label{eq:fgchar}
f_g(w) = f_{q_g}(w) = \sum_{\lam\in\Ghat}
\left(\frac{|G|}{\dim\lam}\right)^{2g-1}\rchi_\lam(w)\;.
\end{equation}
If the sum of the RHS converges with respect to the sup norm on
$G$, then
the value of $f_g(w)$ gives the volume $|q_g ^{-1}(w)|$,
and in particular, we have 
\begin{equation}
\label{eq:volumeformulaori}
|G|^{\rchi(S)-1}|\Hom(\pi_1(S),G)| = 
\sum_{\lam\in\Ghat}(\dim\lam)^{\rchi(S)}\;.
\end{equation}
Conversely, if the RHS of (\ref{eq:volumeformulaori})
is convergent, then the character expansion
(\ref{eq:fgchar}) is uniformly and absolutely convergent
and the volume distribution is a $C^\infty$ class function.
\end{thm}

\begin{rem}
Except for the finite group case, the RHS of
(\ref{eq:volumeformulaori}) never converges for $g=0$ and 
$g=1$. The Weyl dimension formula
gives a lower bound on the genus $g$ for the series to 
converge. For example, if $G=SU(n)$, 
then the series (\ref{eq:volumeformulaori}) 
converges 
 for every $g\ge 2$. 
\end{rem}

We note that the argument works without any modification
for a finite group. 
If $G$ is a finite group, then there is a canonical choice
of the Haar measure: a discrete measure of uniform weight
$1$ for each element. With this choice, the volume $|G|$ is 
the order of the group and the delta function  is the
characteristic function that takes value $1$ at the identity
element and $0$ everywhere else.
There is no difference 
between distributions and differentiable functions
in this case, and there is no problem
of convergence in the character expansion.
 Thus $f_g(w)$ of  (\ref{eq:fgchar})
is a function with a well-defined value for every $w\in G$.
In particular, the expression of $f_g(1)$ for $g\ge 1$ recovers 
Mednykh's formula \cite{Med}
(Theorem~\ref{thm:introori}).
 For a surface of genus $0$,
the formula reduces to the classical formula 
 (\ref{eq:orderformula}) and remains 
true.

For the regularity of the distribution
of $f_g(w)$ at $w=1$, we have the following:

\begin{prop}
\label{prop:regularity}
Let $S$ be a closed orientable surface of genus $g$ 
and $G$ a  compact connected semisimple Lie group 
of dimension $n$. If $g\ge n$, then 
the volume distribution of (\ref{eq:fg}) is
regular at $w=1$.
\end{prop}

\begin{proof}
Since $G$ is semisimple, its Lie algebra
$\Lieg$  satisfies
\begin{equation}
\label{eq:derivedalgebra}
D(\Lieg) = [\Lieg,\Lieg] = \Lieg\;.
\end{equation}
This means every element of $G$ near the identity 
is expressible as a finite product of commutators. 
Let us analyze the map
$$
q_1: G\times G\owns (x,y)\longmapsto xyx^{-1}y^{-1}\in G
$$
near the identity element. Put
$x=e^{tA}$ and $y=e^{sB}$ for some $A,B\in\mathfrak{g}$. 
Then from the Campbell-Hausdorff formula we have
\begin{equation*}
\begin{split}
xyx^{-1}y^{-1} &= e^{tA+sB+\frac{ts}{2}[A,B]+\cdots}
e^{-tA-sB+\frac{ts}{2}[A,B]+\cdots } \\
&=1+ts[A,B] + \text{higher}\;.
\end{split}
\end{equation*}
Therefore, for a fixed $y=e^B$, the tangent line at the
identity of the image of the  one-parameter subgroup
$e^{tA}$ via $q_1$ is given by 
$t[A,B]\in \mathfrak{g}$. 

Now let $(A_1,B_1,\dots,A_g,B_g)\in \Lieg ^{2g}$
be an arbitrary $2g$-tuple of Lie algebra elements. 
The differential of the map $q_g:G^{2g}\rightarrow
G$ at $(1,\dots,1)$ of $G^{2g}$ is 
\begin{equation}
\label{eq:derofqg}
dq_g:\Lieg ^{2g}\owns (A_1,B_1,\dots,A_g,B_g)
\longmapsto [A_1,B_1]+\cdots+[A_g,B_g]\in\Lieg\;.
\end{equation}
If $g\ge n=\dim G$, then $dq_g$ is surjective
because of (\ref{eq:derivedalgebra}). This implies that
the map $q_g$ is regular at $(1,\dots,1)$. 

Let $(a_1,b_1,\dots,a_g,b_g)\in q_g ^{-1}(1)$ be an arbitrary
point of the inverse image of $1$ via $q_g$. Since the exponential 
map is surjective, we can write $a_j=e^{X_j}$ and $b_j=e^{Y_j}$
for $j=1,\dots,g$. Then again by the Campbell-Hausdorff formula,
the differential $dq_g$ at $(a_1,b_1,\dots,a_g,b_g)$ gives
for $(A_1,B_1,\dots,A_g,B_g)\in \Lieg ^{2g}$ the value
$$
\sum_{j=1} ^g ([X_j,B_j]+[A_j,Y_j]+[A_j,B_j]) + \cdots\;,
$$
where we used the fact that $(a_1,b_1,\dots,a_g,b_g)$ is 
in the inverse
image of $1$. 
This is a surjective map onto $\Lieg$.
Therefore, $q_g$ is
locally a fibration in a neighborhood of the identity, and
hence the volume distribution $f_g(w)$ is
regular at $w=1$ by Proposition~\ref{prop:fqregularity}.
\end{proof}

\begin{rem}
\begin{enumerate}
\item 
The semisimplicity condition of  $G$ is necessary.
Indeed, if $G$ has a center $Z$ of positive dimension,
then $Z\times G$ is in $q_1 ^{-1}(1)$, which has a larger
dimension than the expected dimension of the inverse image.
Thus the map $q_g$ has   critical
value at $1$, and  hence $f_g(w)$ may not be 
regular  at $w=1$. 
In particular, the volume distribution is always
singular at the identity for any compact
abelian group of positive dimension. 

\item
Even though $f_g(1)$ is regular, it does not imply
that the series of the RHS of (\ref{eq:volumeformulaori})
 is  convergent. 
\end{enumerate}
\end{rem}

\begin{ex}
Let us consider the case $G=SU(2)$. Since $\Ghat=\mathbb{N}$,
we identify $\lam$ with a positive integer $n$, which is the
dimension of $\lam$. The character expansion (\ref{eq:fgchar}) is
uniformly and absolutely
 convergent for $g\ge 2$, and we have for every closed
orientable surface $S$ of $g\ge 2$,
$$
\frac{|\Hom(\pi_1(S),SU(2))|}{|SU(2)|^{2g-1}} = \zeta(2g-2)\;.
$$
\end{ex}

\bigskip

\section{Representation Varieties of 
Non-orientable Surface Groups}
\label{sect:nonori}

Let us now 
turn our attention to the case of a non-orientable surface $S$.
We recall that every closed non-orientable surface is obtained by
removing $k\ge 1$ disjoint disks from a sphere $S^2$ and gluing
a cross-cap to each hole. The number $k$ is the 
\emph{cross-cap genus} of $S$, and its Euler characteristic is
given by $\rchi(S) = 2-k$. Since the fundamental group of 
a non-orientable surface of cross-cap genus $k$ is given by
\begin{equation}
\label{eq:pi1nonori}
\pi_1(S) = \langle a_1,\dots,a_k\;|\;a_1 ^2 \cdots a_k ^2 \rangle\;,
\end{equation}
we use the presentation map
\begin{equation}
\label{eq:rk}
r_k:G^k\owns (x_1,\dots,x_k)\longmapsto
x_1 ^2\cdots x_k ^2\in G\;.
\end{equation}
Our aim is to compute the volume distribution 
\begin{equation}
\label{eq:hk}
\begin{split}
h_k(w) = f_{r_k}(w) = |r_k ^{-1}(w)| 
&= \int_{G^k} \delta(r_k(x_1,\dots,x_k)\cdot w^{-1})
dx_1\cdots dx_k\\
&=\int_{G^k} \delta(x_1 ^2\cdots x_k ^2 w^{-1})
dx_1\cdots dx_k\;.
\end{split}
\end{equation}
Using the same technique of Section~\ref{sect:ori}, 
it follows from the convolution property of the $\delta$-function
that
$$
h_k = \overset{k\text{-times}}{\overbrace{h_1 *\cdots * h_1}}
\;.
$$

\begin{prop}
\label{prop:h1}
The character expansion of the class distribution $h_1$ 
is given by 
\begin{equation}
\label{eq:h1}
h_1(w) = \int_{G} \delta(x^2 w^{-1})dx = 
\sum_{\lam\in\Ghat} \phi(\lam) \rchi_\lam (w)\;,
\end{equation}
where 
\begin{equation}
\label{eq:fs}
\phi(\lam) = \frac{1}{|G|}\int_G \rchi_\lam(x^2)dx
\end{equation}
is the Frobenius-Schur indicator of \cite{FS}.
\end{prop}

\begin{proof}
Let $h_1(w) = \sum_\lam a_\lam \rchi_\lam(w)$ be
the expansion in terms of  irreducible characters with respect to 
the strong convergence in $\mathcal{D}'(G)$. By definition, we have
\begin{equation*}
\begin{split}
a_\lam 
&= \frac{1}{|G|}\int_G h_1(w)
\overline{\rchi}_\lam(w)dw \\
&=\frac{1}{|G|}\int_G \int_G \delta(x^2w^{-1})
\overline{\rchi}_\lam(w)dwdx\\
&=\frac{1}{|G|}\int_G 
\overline{\rchi}_\lam(x^2)dx\\
&=\frac{1}{|G|}\int_G \rchi_\lam(x^2)dx\;.
\end{split}
\end{equation*}
\end{proof}

\begin{rem}
Since the presentation map $r_1$ has the same domain and the range,
if the space $r_1 ^{-1}(1)$ of involutions of $G$ has positive dimension,
then the map $r_1: G\owns x\mapsto x^2\in G$ does not 
have regular value at $1\in G$. Let $G$ be a compact 
connected Lie group of rank $r$ and $T$  a maximal
torus of $G$. A compact torus  of 
dimension $r$ has exactly $2^r $ involutions. The
conjugacy class of each of these involutions has  positive dimension
unless the involution is central. Thus the groups
$SO(n)$ and $SU(n)$ for $n\ge 3$ and
$Sp(n)$ for $n\ge 2$ all have large
space of involutions with positive dimensions. 
Therefore, the volume $h_1(1)$ is ill-defined
for these groups.
\end{rem}

Taking the $k$-th convolution power of $h_1$, we have

\begin{thm}
\label{thm:hk}
Let $S$ be a closed non-orientable surface of cross-cap genus $k\ge 1$
and $G$ an arbitrary compact Lie group. Then the 
character expansion of the volume distribution is given by
\begin{equation}
\label{eq:hkexpansion}
h_k(w) = |r_k^{-1}(w)| = \sum_{\lam\in\Ghat_1}
\left(\frac{|G|}{\dim\lam}\right)^{k-1} \rchi_\lam(w)
-
\sum_{\lam\in\Ghat_4}
\left(-\;\frac{|G|}{\dim\lam}\right)^{k-1} \rchi_\lam(w)\;.
\end{equation}
If the character sum of the RHS is absolutely convergent
with respect to the sup norm on $G$, 
then the value gives the volume $|r_k ^{-1}(w)|$. In particular,
we have
\begin{equation}
\label{eq:volnonori}
|G|^{\rchi(S)-1} |\Hom(\pi_1(S),G)| =  
 \sum_{\lam\in\Ghat_1}
(\dim\lam)^{\rchi(S)} 
+
\sum_{\lam\in\Ghat_4}
(-\dim\lam)^{\rchi(S)}\;.
\end{equation}
Conversely, if the RHS of (\ref{eq:volnonori}) is absolutely
convergent, then (\ref{eq:hkexpansion}) is uniformly and
absolutely convergent and the volume distribution 
$h_k$ is a $C^\infty$
class function.
\end{thm}

\begin{proof}
All we need is to notice the value of the
Frobenius-Schur indicator (\ref{eq:FSindicator}). 
\end{proof}

\begin{ex}
\label{ex:torus}
Consider an $n$-dimensional torus $G=T^n$. 
Since the only real irreducible representation is the trivial
representation, $\Ghat_1=\{1\}$. 
Let $S^1 = \{e^{i\theta}\;|\; 0\le \theta < 2\pi\}$
and choose an invariant measure 
$d\theta$ on it. The $\delta$-function
on the torus is the product of $\delta$-functions on $S^1$. 
Thus we have
\begin{equation*}
\begin{split}
h_1(w)&=\int_{T^n} \delta(2\theta_1 -w_1)\cdots\delta(2\theta_n-w_n)
d\theta_1\cdots d\theta_n\\
&=\left(\int_0 ^{2\pi} \delta(2\theta)d\theta\right)^n\\
&=\left(\frac{1}{2}\int_0 ^{4\pi} \delta(\theta)d\theta\right)^n\\
&=1\;.
\end{split}
\end{equation*}
Therefore, $h_k(w)=|T^n|^{k-1}$, and hence 
$|\Hom(\pi_1(S),T^n)| = |T^n|^{1-\rchi(S)}$.  
The inverse image $r_k ^{-1}(w)$ consists of 
$(x_1,\dots,x_k)\in (T^n)^k$
such that 
\begin{equation}
\label{eq:torus}
x_k ^2 = x_{k-1} ^{-2}\cdots x_{1} ^{-2}w\;.
\end{equation}
Thus we can choose arbitrary $(x_1,\dots,x_{k-1})\in (T^n)^{k-1}$,
and for this choice, there are still $2^n$ solutions of (\ref{eq:torus}). 
Each piece has the volume equal to $|T^n|^{k-1}/2^n$. Therefore, the total 
volume is $|T^n|^{k-1}$, in agreement with our computation. 
\end{ex}

\begin{ex}
For $G=SU(2)$, we have a more interesting relation. We note that
$\Ghat_1$ consists of odd integers and $\Ghat_4$ even integers. 
Therefore, the series (\ref{eq:volnonori}) is 
absolutely convergent for $k\ge 4$. Thus for a closed non-orientable
surface of cross-cap genus $k\ge 4$, we have
\begin{equation*}
\begin{split}
\frac{|\Hom(\pi_1(S),SU(2))|}{|SU(2)|^{k-1}}
&=\sum_{n=1} ^\infty (2n-1)^{2-k} + (-1)^{2-k}\sum_{n=1} ^\infty
(2n)^{2-k}\\
&=
\left\{
\begin{matrix}
\zeta(k-2)&& k \text{ is even}\\
\left(1-\frac{1}{2^{k-3}}\right) \zeta(k-2) && k \text{ is odd}
\end{matrix}
\right. \;.
\end{split}
\end{equation*}
Although the series is not absolutely convergent, we expect
that the value for $k=3$ is $\log(2)$. 
\end{ex}

In the same way as in Proposition~\ref{prop:regularity},
we have a general regularity condition of
the volume distribution $h_k(w)$ for a 
compact semisimple Lie group.

\begin{prop}
\label{prop:hkregularity}
Let $G$ be a compact connected semisimple Lie group of 
dimension $n$ and $S$ a closed non-orientable surface of
cross-cap genus $k\ge 1$. Define 
\begin{equation*}
m=\left[\frac{k-1}{2}\right]
=\begin{cases}
(k-1)/2\qquad k\; \text{ is odd}\\
(k-2)/2\qquad k\; \text{ is even.}
\end{cases}
\end{equation*}
Then the volume distribution $h_k(w)$ is regular at $w=1$
if $m\ge n$.
\end{prop}

\begin{proof}
The argument is based on another presentation of the 
fundamental group
 containing
$m$ commutator products
$$
\pi_1(S) = \left\langle a_1,b_1,\dots,a_m,b_m,c\;\left|\;
[a_1,b_1]\cdots [a_m,b_m]\; c^2\right\rangle\right.
$$
for odd $k$ and
$$
\pi_1(S) = \left\langle a_1,b_1,\dots,a_m,b_m,c_1,c_2\;\left|\;
[a_1,b_1]\cdots [a_m,b_m]\; c_1 ^2 c_2 ^2\right\rangle\right.
$$
for even $k$.
First we notice that the substitution
\begin{equation}
\label{eq:substitution}
\begin{cases}
\alpha =abc\\
\beta = c^{-1}b^{-1}a^{-1}c^{-1}a^{-1}c\\
\gamma = c^{-1}ac^2\;,
\end{cases}
\end{equation}
yields the equality
\begin{equation}
\label{eq:m1case}
[a,b]\;c^2 = \alpha^2 \beta^2 \gamma^2\;.
\end{equation}
The inverse transformation of (\ref{eq:substitution})
is given by 
\begin{equation}
\label{eq:inverse}
\begin{cases}
a=\alpha \beta \gamma \beta^{-1} \alpha^{-1} \gamma^{-1}
\beta^{-1} \alpha^{-1}\\
b=\alpha \beta \gamma \alpha \beta \gamma^{-1} \beta^{-1}
\gamma^{-1} \beta^{-1} \alpha^{-1}\\
c=\alpha\beta\gamma\;.
\end{cases}
\end{equation}
Therefore, the map
$$
h:G^3\owns (x,y,z)\longmapsto 
(xyz,z^{-1}y^{-1}x^{-1}z^{-1}x^{-1}z,
z^{-1}xz^2)=(u,v,w)\in G^3
$$
defined by (\ref{eq:substitution}) is an analytic automorphism of
the manifold $G^3$. Since both $h$ and $h^{-1}$ are
given by words,  from Lemma~\ref{lem:volumeforminv} we conclude
\begin{equation}
\label{eq:homegaisomega}
h^* \omega^3 = \omega^3\;.
\end{equation}
Indeed, the integer matrix $N$ representing
the differential $dh$ (see proof of Lemma~\ref{lem:volumeforminv})
is 
\begin{equation*}
N=
\begin{bmatrix}
1&1&1\\
-2&-1&-1\\
1&0&1
\end{bmatrix}
\end{equation*}
with $\det N = 1$. 
Let
$$
s_3:G^3\owns (x,y,z)\longmapsto [x,y]\; z^2\in G
$$
be the presentation map with respect to the relator $[a,b]\;c^2$, 
and 
$$
r_3:G^3\owns (u,v,w)\longmapsto u^2 v^2 w^2\in G
$$
the presentation map of $\alpha^2 \beta^2 \gamma^2$. 
Then (\ref{eq:m1case}) implies that the following diagram
commutes:
\begin{equation*}
\begin{CD}
G^3 @>{s_3}>> G\\
@V{h}V{\wr}V @|\\
G^3 @>>{r_3}> G
\end{CD}\;.
\end{equation*}
Together with the invariance of the volume form
(\ref{eq:homegaisomega}), we obtain
$$
f_{s_3}(w) = f_{r_3}(w) = h_3(w)
$$
as a distribution on $G$. 

We can apply the above procedure based on the transformation
(\ref{eq:substitution}) repeatedly to change the relator
$$
[a_1,b_1]\cdots [a_m,b_m] \; c^2
\qquad \text{to}\qquad 
\alpha_1 ^2\beta_1 ^2\cdots \alpha_m ^2 \beta_m ^2 \gamma^2\;,
$$
and 
$$
[a_1,b_1]\cdots [a_m,b_m] \; c_1 ^2 c_2 ^2
\qquad \text{to}\qquad
\alpha_1 ^2\beta_1 ^2\cdots \alpha_m ^2 \beta_m ^2 
\gamma_1 ^2\gamma_2 ^2\;. 
$$At each step the volume form
on $G^k$ is kept invariant, and we have the equality of the
volume distributions
\begin{equation}
\label{eq:fish}
f_{s_k}(w) = h_k (w)\;,
\end{equation}
where $s_k:G^k\longrightarrow G$ is the presentation 
map associated with the relator containing $m$ commutators. 
The regularity of $f_{s_k}(w)$ at $w=1\in G$
can be established in the same
way as  in Proposition~\ref{prop:regularity}.
The  regularity of $h_k(w)$ follows from (\ref{eq:fish}).
\end{proof}
\bigskip

\section{Harmonic Analysis and 
Irreducible Representations of a Compact Lie Group}
\label{sect:harmonic}

In this section we briefly review the necessary
accounts from the character theory of compact Lie groups
and non-commutative harmonic analysis that are used in the
main part of this paper. 
The important results we need are Schur's orthogonality
relations of irreducible characters in terms of the convolution 
product, definition of the $\delta$-function on a compact Lie group,
and expansion of a class distribution in terms of irreducible 
characters. We refer to \cite{BD, Gelfand, V, V2, Warner} for
more detail.

Let $G$ be a compact Lie group.  We choose a left and right invariant,
but not necessarily normalized, 
measure $dx$ (Haar measure) on $G$ that satisfies 
\begin{equation}
\label{eq:Haar}
\int_G g(x)dx = \int_G g(yx)dx =\int_G g(xy)dx =\int_G g(x^{-1})dx\;,
\end{equation}
where $g(x)$ is a smooth function on $G$ and $y\in G$. 
The space of \emph{distributions} $\mathcal{D}'(G)$ is the
topological linear space
consisting of all continuous linear functions on $C^\infty (G)$
with respect to the Fr\'echet topology of $C^\infty (G)$
\cite{Gelfand, Warner}. 
The left regular representation of $G$ on $\mathcal{D}'(G)$
is defined by
\begin{equation}
\label{eq:leftreg}
\rho_L(w):\mathcal{D}'(G)\owns f(x)
\longmapsto f(w^{-1}x)\in \mathcal{D}'(G)\;,
\end{equation}
where the distribution $f(w^{-1}x)$ is characterized  by
\begin{equation}
\label{eq:fwx}
\int_G f(w^{-1} x) g(x)dx = 
\int_G f(x) g(wx)dx
\end{equation}
for every $g(x)\in C^\infty(G)$
using the invariance (\ref{eq:Haar}) of $dx$. Similarly,
the right regular representation is defined by
\begin{equation}
\label{eq:rightreg}
\rho_R(w):\mathcal{D}'(G)\owns f(x)
\longmapsto f(xw)\in \mathcal{D}'(G)\;.
\end{equation}
Reflecting the group structure of the compact space $G$,
$\mathcal{D}'(G)$ has the structure of an associative
algebra over $\mathbb{C}$ whose multiplication is the 
\emph{convolution product}
of distributions $f_1$ and $f_2$ defined by
\begin{equation}
\label{eq:conv}
\big(f_1*f_2\big)(x) = \int_G f_1(xw^{-1}) f_2(w) dw\;.
\end{equation}
A distribution $f(x)\in \mathcal{D}'(G)$ is a
\emph{class distribution} if  
\begin{equation}
\label{eq:classfunc}
\rho_L(w^{-1})\cdot f(x) = \rho_R(w)\cdot f(x)\;,
\end{equation}
or more conveniently, when
$f(w^{-1}xw) = f(x)$ as a distribution for every 
$w\in G$.
Class distributions form the center of the 
convolution algebra $\mathcal{D}'(G)$ because of the
invariance of the Haar measure. Indeed, for every
smooth function $g(x)\in C^\infty(G)$, the condition
$$
\big(f*g\big)(x)=\int_G f(xw^{-1})g(w)dw = 
\int_G g(xu^{-1}) f (u)du = \big(g*f\big)(x)
$$
for $f(x)\in\mathcal{D}'(G)$ is equivalent to (\ref{eq:classfunc})
because
$$
\int_G g(xu^{-1}) f (u)du =
\int_G f(w^{-1}x) g(w)dw\;.
$$
The Dirac $\delta$-function on $G$ is defined by 
\begin{equation}
\label{eq:delta}
\int_G \delta(x)g(x)dx = g(1)\;,\qquad g(x)\in C^\infty(G)\;.
\end{equation}
This definition, together with the left and right invariance of
the Haar measure $dx$, implies the following formulas
for an arbitrary $y\in G$:
\begin{equation}
\label{eq:formulasofdelta}
\begin{split}
\int_G \delta(xy^{-1})f(x)dx &= f(y)\;, \qquad f(x)\in \mathcal{D}'(G)\\
\delta(yxy^{-1})&=\delta(x)=\delta(x^{-1})\\
\delta(xy)&=\delta(yx)\;.
\end{split}
\end{equation}
Therefore, the $\delta$-function is a class distribution
and serves as the identity element of the convolution 
algebra $(\mathcal{D}'(G),*)$. In particular, we have
$\delta *\delta = \delta$.

The Hermitian inner product of $L^2$ functions
$f_1$ and $f_2$ on $G$ is given by
\begin{equation}
\label{eq:herm}
\langle f_1,f_2\rangle = \frac{1}{|G|}
\int_G f_1(x) \overline{f_2(x)}dx\;.
\end{equation}
Peter-Weyl theory provides a Hilbert basis for $L^2(G)$. 
The $L^2$ class functions on $G$ form
a Hilbert space with a Hilbert basis
\begin{equation}
\label{eq:hilbertbasis}
\{\rchi_\lam\;|\;\lam\in\Ghat\}\;,
\end{equation}
where $\Ghat$ is the set of equivalence classes of complex
irreducible representations of $G$ and
$\rchi_\lam$ the character of $\lam\in\Ghat$. 
The irreducible characters 
satisfy \emph{Schur's orthogonality relations}
\begin{equation}
\label{eq:ortho}
\langle \rchi_\lam , \rchi_\mu\rangle =
 \delta_{\lam\mu} 
\end{equation}
and
\begin{equation}
\label{eq:convortho}
\rchi_\lam * \rchi_\mu =
 \frac{|G|}{\dim\lam}\; \delta_{\lam\mu}
\rchi_\lam\; ,
\end{equation}
where $\dim\lam$ is the dimension of the representation
$\lam\in\Ghat$. 

Since $G$ is compact, every finite dimensional representation
has a natural $G$-invariant Hermitian inner product so that 
the representation is unitary. In particular, we have
\begin{equation}
\label{eq:estimateofchar}
|\rchi_\lam(x)|\le |\rchi_\lam(1)| = \dim \lam
\qquad \text{ for every } x\in G\;.
\end{equation}
Therefore, if a sequence $a_\lam\in\mathbb{C}$ satisfies 
the convergence condition
\begin{equation}
\label{eq:conditiona}
\sup_{x\in G}\left|\sum_{\lam\in\Ghat} a_\lam \rchi_\lam (x)\right|\le 
\sum_{\lam\in\Ghat} |a_\lam|\dim \lam  <\infty\;,
\end{equation}
then the series $\sum_\lam a_\lam \rchi_\lam$ is uniformly and 
absolutely convergent and defines a 
$C^\infty$ class function on $G$. 
Note that every $L^2$ class function on $G$
can be expanded into a series 
of the form $\sum_\lam a_\lam \rchi_\lam$ with the convergence
condition $\sum_\lam |a_\lam|^2<\infty$.
(If $G$ is a compact connected abelian group, i.e., a 
compact torus, then the expansion in terms of 
irreducible characters is the Fourier 
series expansion.) 
In the same way, a class distributions on $G$ is
expanded into a series in $\rchi_\lam$ 
that is not necessarily absolutely convergent. For example,
the $\delta$-function has an expansion
\begin{equation}
\label{eq:deltaexpansion}
\delta(x) = \frac{1}{|G|}\sum_{\lam\in\Ghat}
\dim\lam \cdot \rchi_\lam(x)\;.
\end{equation}
This infinite sum is convergent with respect to the 
strong topology of the topological linear space
$\mathcal{D}'(G)$, and gives an alternative definition of the delta
function.
We note that the convolution relation $\delta*\delta=\delta$
is equivalent to Schur's orthogonality relation
 (\ref{eq:convortho}).

The following  family of 
class distributions parametrised by $x\in G$ are used 
in Section~\ref{sect:ori}:
\begin{equation}
\label{eq:eta}
\eta_x(w) = \int_G \delta(yxy^{-1}w^{-1})dy\;.
\end{equation}
It is the distribution corresponding to  the linear functional
$$
g(w)
\longmapsto
\int_G g(yxy^{-1})dy
$$
for every smooth function $g(w)$. 
It is a class distribution because for every $u\in G$,
$$
\int_G g(yxy^{-1})dy
= \int_G g(uyxy^{-1}u^{-1})d(uy)
=\int_G g(u(yxy^{-1})u^{-1})dy\;.
$$
The character
expansion of $\eta_x(w)$ as a distribution in $w$ has 
the following simple form:
\begin{equation}
\label{eq:expansionofeta}
\eta_x(w) = \sum_{\lam\in\Ghat}
\overline{\rchi}_\lam(x)\rchi_\lam(w)\;.
\end{equation}

Distributions cannot be evaluated in a meaningful way in general.
For example, $\delta(1)$ is not defined. Similarly, 
restriction of a distribution defined on $G$ to a closed
subset  is generally meaningless. However,  a distribution can 
always be localized to any open subset of $G$. Let
$U\in G$ be an open subset and  $\mathcal{D}(U)$
the Fr\'echet space of $C^\infty$ functions on $U$ with compact
support. The \emph{localization} of a distribution on $G$
to $U$, 
\begin{equation}
\label{eq:localization}
\mathcal{D}'(G)\owns f\longmapsto
f|_U\in \mathcal{D}'(U)\;,
\end{equation}
is defined by
$$
\int_U \big(f|_U\big)(x)g(x)dx = \int_G f(x)\tilde{g}(x) dx
$$
for every $g(x)\in \mathcal{D}(U)$, where 
$\tilde{g}(x)$ is the extension 
of $g(x)$ as a $C^\infty$ function
on $G$ satisfying
\begin{equation*}
\tilde{g}(x)
=\left\{
\begin{matrix}
g(x) & x\in U\\
0& x\notin U
\end{matrix}
\right.\;.
\end{equation*}
Since $C^\infty (U)\subset \mathcal{D}'(U)$, we can 
compare $f|_U$ with a smooth function on $U$. 
A distribution $f\in \mathcal{D}'(G)$ is said to be
\emph{regular} at $w\in G$ if there is an open 
neighborhood $U\owns w$ such that the localization
$f|_U$ of $f$ is a $C^\infty$ function on $U$. 
If $f$ is regular at $w\in G$, then the value 
$f(w)\in \mathbb{C}$ of $f$ at $w$  is well-defined.

\bigskip

\section{Representation Varieties and Cell-Decompositions
of a Surface}
\label{sect:cell}

We have observed that for every closed surface
$S$ and a compact Lie groups $G$, the quantity
$|G|^{\rchi(S)-1}|\Hom(\pi_1(S),G)|$ is a natural
object that does not depend on the choice of an
invariant measure on $G$. In a recent work
\cite{MY, MY2}, it has been noted that the 
same quantity for a finite group
appears in the graphical expansion of a
certain integral over the group algebra of the finite group.
In this section we present a \emph{heuristic} 
argument that relates the volume of the representation 
variety of a surface group to each term of
the Feynman diagram expansion 
of a hypothetical ``integral'' over the group algebra of
compact Lie groups. The argument is based on an
alternative definition of the volume of the representation 
variety of a surface group in terms of a cell-decomposition
of the surface.

The $1$-skelton of a cell-decomposition of a closed connected 
surface is a \emph{graph} drawn on the surface. Such graphs
are known as  \emph{maps}, \emph{fatgraphs}, \emph{ribbon graphs},
\emph{dessins d'enphants}, or \emph{M\"obius graphs}, 
depending on the context of study.
Fatgraphs and ribbon graphs
are used effectively in the study of topological 
properties of moduli spaces of pointed Riemann
surfaces and related topics (see for example,
\cite{BIZ, Harer1986, Harer-Zagier, Kontsevich,
O, OP1, OP2, Penner, W1991a}). 
Dessins d'enphants appear in the context of algebraic curves
defined over the field of algebraic numbers
(cf.~\cite{Belyi, Grothendieck, MP1998, MP2001,
Schneps, SL}). Maps are drawn on
 a topological surface 
and are studied from the point of view of 
map coloring theorems (see for example, \cite{GT, Ringel}),
where non-orientable surfaces are also considered.
Graphs on non-orientable surfaces play an important
role in the study of moduli spaces of real algebraic 
curves \cite{GHJ}. In this section we use the notion of 
\emph{M\"obius graphs} 
following \cite{MW, MY2} to emphasize that our
graphs are on \emph{unoriented} surfaces, which
may or may not be orientable.

A ribbon graph is a graph with a cyclic order given at
each vertex to the set of half-deges incident to it.
Equivalently, it is a graph drawn on a closed oriented 
surface \cite{MP1998}. Vertices with cyclic orders
 are connected by ribbon like
edges to form a ribbon graph \cite{Kontsevich}. 
In order to deal with a graph on a non-orientable surface,
it is necessary to allow ribbons with a \emph{twist}. 
We can use the group $\mathfrak{S}_2 = \{\pm 1\}$ to 
indicate a twist of an edge with the value $-1$ and no twist 
with $1$.
We are thus led to 
 consider a ribbon graph with $\mathfrak{S}_2$-color
at every edge. 
A \emph{vertex flip} operation of an edge-colored ribbon graph
is a move at a vertex
$V$ of a graph that reverses the cyclic order at $V$ and 
interchanges the color of the edge that is incident to $V$. 
If an edge makes a loop and incident to $V$ twice, then 
its color is changed twice and hence preserved. 
Two $\mathfrak{S}_2$-colored ribbon graphs are
\emph{equivalent} if one is brought to the other by a sequence of
vertex flip operations. 
A \emph{M\"obius graph} is the equivalence class of a
$\mathfrak{S}_2$-colored ribbon graph \cite{MY2}.

Let $\Gamma_c$ be a  connected
$\mathfrak{S}_2$-colored ribbon graph,
and $G$ a compact Lie group. For every vertex $V$ of 
$\Gamma_c$, we label the half-edges incident to $V$ as
$(V1,V2,\dots,Vn)$ according to the cyclic order, where $n$ is the
valence (or degree) of the vertex. 
We call the expression 
\begin{equation}
\label{eq:vertexcontribution}
\Phi(V)=
\delta(x_{V1}\cdots x_{Vn}) \in \mathcal{D}'(G^n)
\end{equation}
the \emph{vertex contribution}. 
It is invariant under the cyclic permutation of the half-edges.
Suppose an edge $E$ is incident
to  vertices $V$ and $W$ of $\Gamma_c$, and consists of two 
half-edges $Vi$ and $Wj$. The \emph{propagator} of the edge $E$
is the expression
\begin{equation}
\label{eq:propagator}
\Phi(E)=
\delta(x_{Vi} ^{}  x_{WJ} ^c)\in \mathcal{D}'(G^2)\;,
\end{equation}
where $c\in \mathfrak{S}_2=\{\pm 1\}$ is the color of $E$.
A few remarks are necessary for the definition of
a vertex contribution and a propagator as distributions. 
For every $n\ge 2$, the multiplication
$$
m_n: G^n \owns (x_1,\dots,x_n)\longmapsto
x_1\cdots x_n\in G
$$
is a real analytic map and 
$$
\psi_n:G^n \owns (x_1,\dots,x_{n-1},x_n)
\overset{\sim}{\longmapsto}
\big( (x_1\dots,x_{n-1}),(x_1\cdots x_{n-1})\cdot x_n\big)\in
G^{n-1}\times G
$$
is a real analytic isomorphism. 
The Jacobian $J(x_1,\dots,x_{n-1},y)$
of the isomorphism $\psi_n$
gives the push-forward measure
$$
(\psi_n)_*(dx_1\cdots dx_n)=
J(x_1,\dots,x_{n-1},y)dx_1\cdots dx_{n-1} dy\;.
$$
The vertex contribution (\ref{eq:vertexcontribution}) is 
a distribution defined by the continuous linear functional
\begin{multline*}
C^\infty(G^n)\owns f(x_1,\dots,x_n)
\longmapsto \\
\int_{G^{n-1}} 
J(x_1,\dots,x_{n-1},1)
f\big( x_1,\dots,x_{n-1},(x_1\cdots x_{n-1})^{-1}
\big) dx_1\cdots dx_{n-1}
\in\mathbb{C}\;.
\end{multline*}
The propagator as a distribution is defined similarly. 
We note that the product of all vertex contributions and
the product of all propagators 
$$
\prod_{V\in \Gamma_c} \Phi(V) 
\qquad\qquad
\text{and}
\qquad\qquad
\prod_{E\in \Gamma_c} \Phi(E)
$$
are well-defined distributions in 
$\mathcal{D}'(G^{2e})$,
where $e=e(\Gamma_c)$ is the number of edges of 
the graph.

Now we come to the heuristic part of the argument. 
Let us ``define''
$$
G_{\Gamma_c} = 
\int_{G^{2e(\Gamma_c)}}
\prod_{V\in \Gamma_c} \Phi(V)
\prod_{E\in \Gamma_c} \Phi(E)
dx\;,
$$
where $dx$ is the product measure on 
$G^{2e(\Gamma_c)}$. Although these products 
are individually well-defined distributions on 
the space $G^{2e(\Gamma_c)}$, their product 
does not make sense in general.
For now, let us suppose that $G_{\Gamma_c}$ is 
legally defined. Then the first thing we can show is that
\emph{it is invariant under vertex flip operations}. 
To see this, choose a vertex $V$ at which we perform the vertex
flip. Suppose it is incident to half-edges $V1,\dots,Vn$ in this 
cyclic order. We denote by $V'i$ the other half-dege
connected to $Vi$. Of course we do not mean they are 
all incident to another vertex labeled by $V'$. This is just
a notation for the other half. Then the vertex flip invariance
follows from
\begin{equation*}
\begin{split}
\delta(x_{V'1} ^{-c}\cdots x_{V'n} ^{-c})
&=\int_{G^n}
\delta(x_{V1}\cdots x_{Vn})\delta(x_{V1}x_{V'1} ^c)
\cdots \delta(x_{Vn}x_{V'n} ^c)
dx_{V1}\cdots dx_{Vn}\\
&= \int_{G^n}
\delta(x_{Vn}\cdots x_{V1})\delta(x_{V1}  x_{V'1} ^{-c})
\cdots \delta(x_{Vn} x_{V'n} ^{-c})
dx_{V1}\cdots dx_{Vn}\\
&=\delta(x_{V'n} ^c \cdots x_{V'1} ^{c})\;,
\end{split}
\end{equation*}
where the first and the last delta functions are the same because 
of (\ref{eq:formulasofdelta}).
Notice that the vertex flip invariance of $G_{\Gamma_c}$ 
implies that this quantity is associated to a M\"obius
graph. So we  define
\begin{equation}
\label{eq:GGamma}
G_{\Gamma} = 
\int_{G^{2e(\Gamma)}}
\prod_{V\in \Gamma} \Phi(V)
\prod_{E\in \Gamma} \Phi(E)
dx
\end{equation}
for every M\"obius graph $\Gamma$. 

The
second property we can show, still assuming 
(\ref{eq:GGamma}) being well-defined, is that
\emph{it is invariant under an edge contraction}.
As explained in \cite{MP1998,MY2}, an edge contraction
of a graph removes an edge that is incident to two
distinct vertices and put these vertices together. 
If the edge is not twisted, then this is exactly the same
edge contraction of \cite{MP1998}. If the edge is
twisted, then we first apply a vertex flip operation to
one of the vertices to untwist the edge, and then proceed with
the usual edge contraction. To see the invariance,
suppose $E$ is the edge to be removed, which
is incident to two vertices $V\ne W$. Without loss of generality,
we assume that $E$ is not twisted, as mentioned above.
Let us denote by $\delta(xy)$ the propagator.
Using cyclic permutations of vertex contributions,
we can assume that $\delta(vx)$ and $\delta(yw)$
are the vertex contributions at $V$ and $W$, respectively,
where $v$ and $w$ are products of group elements. 
The invariance under the edge contraction then follows from
\begin{equation}
\label{eq:edgecontraction}
\int_{G^2}\delta(vx)\delta(xy)\delta(yw)dxdy
= \int_G \delta(vy^{-1})\delta(yw)dy = \delta(vw)\;.
\end{equation}
Note that $\delta(vw)$ is the vertex contribution of the new
vertex obtained by joining $V$ and $W$ together.

Recall the fact that every connected M\"obius graph
$\Gamma$ determines a unique closed surface $S_\Gamma$
and its cell-decomposition. Let $v(\Gamma)$,
$e(\Gamma)$, and $f(\Gamma)$ denote the number of
$0$, $1$, and $2$-cells, respectively. The Euler characteristic
of the surface is given by
$$
\rchi(S_\Gamma) = v(\Gamma)-e(\Gamma)+f(\Gamma)\;.
$$
If $\Gamma$ is orientable, then $S_\Gamma$ is an 
orientable surface of genus $g$ such that $\rchi(S_\Gamma)
=2-2g$,  and if $\Gamma$ is non-orientable, then
$S_\Gamma$ is a non-orientable surface of cross-cap genus
$k$ with $\rchi(S_\Gamma) = 2-k$. Notice that the
edge contraction preserves the Euler characteristic and
the number $f(\Gamma)$ of $2$-cells. Using the same
argument of \cite{MY2}, we can calculate $G_\Gamma$
for every graph. If the graph is orientable of genus
$g\ge 1$, then
\begin{equation}
\label{eq:GGammaori}
G_\Gamma = 
|G|^{f(\Gamma)-1}
\int_{G^{2g}}\delta\big([x_1,y_1]\cdots 
[x_g,y_g]\big)dx_1dy_1\cdots dx_gdy_g\;,
\end{equation}
and if it is non-orientable of cross-cap genus $k$, then
\begin{equation}
\label{eq:GGammanonori}
G_\Gamma = 
|G|^{f(\Gamma)-1}
\int_{G^{k}}\delta\big(x_1 ^2
\cdots x_k ^2\big)dx_1\cdots dx_k\;.
\end{equation}
Regardless the orientation of the graph, we have a general formula
\begin{equation}
\label{eq:GGammageneral}
G_\Gamma = |G|^{f(\Gamma)-1}
\big|\Hom\big(\pi_1(S_\Gamma),G\big)\big|
\end{equation}
if the expressions (\ref{eq:GGammaori}) and (\ref{eq:GGammanonori})
are well-defined. 
In the sense of (\ref{eq:GGammageneral}), we have the
alternative  ``definition''  (\ref{eq:GGamma})
of the volume of the representation variety.

It is clear from the invariance under edge
contraction that $G_\Gamma$ of
(\ref{eq:GGamma}) is always ill-defined 
if $\rchi(S_\Gamma)=2$. As we have seen in Section~\ref{sect:nonori},
the integral (\ref{eq:GGammanonori}) is ill-defined  in general
for graphs
with $\rchi(S_\Gamma)=1$.
For an orientable graph of genus $g\ge 2$, the integral
(\ref{eq:GGammaori}) makes sense only when the group $G$
is semisimple (see Section~\ref{sect:ori}).

Everything we discussed in this section makes sense for 
a finite group $G$, and we have the Feynman diagram expansion formula
\begin{equation}
\label{eq:Feynman}
\log \int_{\{x\in \mathbb{R}[G]\;|\;x^*=x\}}
e^{-\frac{1}{4}\langle x^2\rangle}
e^{\sum_j \frac{t_j}{2j}\langle x^j\rangle}d\mu(x)
=\sum_{\substack{
\Gamma \; \text{ connected}\\
\text{M\"obius graph}}}
\frac{1}{|\Aut_M(\Gamma)|} G_\Gamma
\prod_j t_j ^{v_j(\Gamma)}
\end{equation}
established in \cite{MY,MY2},
where $\langle \;\rangle = \frac{1}{|G|}\rchi_\reg$ is the normalized
trace of the finite-dimensional von Neumann algebra $\mathbb{R}[G]$,
which is the real group algebra of $G$, and $v_j(\Gamma)$ is the number
of $j$-valent vertices of $\Gamma$. We note that if we use
the regular representation itself instead of  the normalized trace in the
integral, then $G_\Gamma$ of the graphical expansion is
replaced by the familiar $|G|^{\rchi(S_\Gamma)} \big|\Hom\big(
\pi_1(S_\Gamma),G\big)\big|$. 
We refer to the above mentioned articles for 
more explanation of the expansion formula.
The point we wish to make here is that the RHS of 
(\ref{eq:Feynman}) \emph{tends to be}
 well-defined for a compact semisimple Lie
group $G$ when the Euler characteristic $\rchi(S_\Gamma)$ is
negative. 
Then how should we consider the integral of the LHS for
a compact Lie group? Is the ill-defined nature of the
infinite dimensional integral concentrated on the non-hyperbolic surfaces?
Namely, if we remove the contributions of
the surfaces of non-negative Euler characteristics
from the integral over the infinite-dimensional space, would then
the integral become well-defined?
If the integral over the real group algebra  $L^2 _\mathbb{R}(G)$
or even $\mathcal{D}' _\mathbb{R}(G)$ could make sense in
some way, then Peter-Weyl theory gives a decomposition 
of the integral into an infinite sum of matrix integrals, and the volume 
formulas  (\ref{eq:volumeformulaori}) and 
(\ref{eq:volnonori}) should follow from the comparison of terms 
corresponding to
 orientable and non-orientable surfaces in the Feynman 
diagram expansion of the integral,
as it was the case for a finite group \cite{MY,MY2}.

The alternative definition (\ref{eq:GGamma}) suggests that 
this subject is related to the study of tensor categories and
modular functors \cite{BK, CP}.

\bigskip


\providecommand{\bysame}{\leavevmode\hbox to3em{\hrulefill}\thinspace}

\bibliographystyle{amsplain}

\end{document}